# Preliminaries on CAT (0) Spaces and Fixed Points of a Class of Iterative Schemes

Manuel De la Sen

Institute of Research and Development of Processes, Leioa 48940, Spain
manuel.delasen@ehu.eus

**Abstract**: This paper gives some relating results for various concepts of convexity in metric spaces $(X,d)$ such as $1$-convexity, $p$-convexity, midpoint convexity, convex structure, uniform convexity and near-uniform convexity, Busemann curvature and its relation to convexity. Some properties of uniform convexity and near uniform convexity of geodesic metric spaces are related to the mapping $Q_t : X \times X \to X$ parameterized by the scalar $t \in [0,1]$, defined by $Q_t(Sx, Ty) = tSx \oplus (1-t)Ty$ in the metric space $(X, d)$ where $S, T : X \to X$ are Lipschitz continuous.

## 1. Introduction

This paper provides some results convexity in metric spaces $(X,d)$ such as $1$-convexity, $p$-convexity, midpoint convexity, convex structure, uniform convexity and near-uniform convexity, Busemann curvature and its relation to convexity. Certain properties of uniform convexity and near uniform convexity of geodesic metric CAT (0)-spaces considers the mapping $Q_t : X \times X \to X$ parameterized by the scalar $t \in [0,1]$, defined by $Q_t(Sx, Ty) = tSx \oplus (1-t)Ty$ in metric space $(X, d)$ where $S, T : X \to X$ are Lipschitz continuous while not necessarily contractive mappings, i.e. the Lipschitz constants are not necessarily less than unity, which is a $CAT(0)$ space. Basic notatiobn used is the following one: $\mathbf{Z}_{0+} = \{z : z \geq 0\}$, $\mathbf{Z}_+ = \{z : z > 0\}$, $\mathbf{R}_{0+} = \{z : z \geq 0\}$, $\mathbf{R}_+ = \{z : z > 0\}$, $c\ell A$ is the closure of the set $A$, $con\{x_i\}_{i \in I}$ is the closure of the convex hull of the family $\{x_i\}_{i \in I}$, $Fix(T)$ denotes the set of fixed points of a mapping $T : X \to X$.

## 2. Preliminaries

Let $(X, d)$ be a complete metric space being a geodesic space.

**Definition 2.1**($p$- convexity) [1, 3]. Suppose a metric space $(X,d)$ which admits midpoints (or which has midpoints or which is midpoint convex). Then, $(X,d)$ is said to be $p$- convex for some $p \in [1, \infty]$ if for each $x, y, z \in X$, and each midpoint $m(x,y) \in X$ of $x$ and $y$:

$$d(m(x,y), z) \leq \left(\frac{1}{2}\right)^{1/p} \left(d^p(x,z) + d^p(y,z)\right)^{1/p} \tag{1}$$

For the case $p = \infty$, the right-hand-side of (1) is defined as a limit leading to $d(m(x,y), z) \leq max(m(x,y), z)$. If $(X,d)$ is $\infty$-convex it is equivalently said to be ball convex while if it is $1$-convex it is equivalently said to be distance convex, [3]. $(X,d)$ is said to be strictly $p$- convex for $p \in (1, \infty]$ if the inequality is strict for $x \neq y$ and strictly 1-convex if the inequality is strict for $p = 1$ if $d(x,y) > |d(x,z) - d(y,z)|$, [1].

**Assertion 2.2**. If a metric space $(X,d)$ is midpoint convex then it is 1-convex. □

**Assertion 2.3**. If $(X,d)$ be $p$-convex for some $p \in [1, \infty)$ then for any $x, y, z \in X$ and each midpoint $m(x,y) \in X$ of $x$ and $y$:

$$d(m(x,y), z) \leq \left(\frac{1}{2}\right)^{1/p} (d(x,z) + d(y,z)) \leq \frac{1}{2}(d(x,z) + d(y,z)) \tag{2}$$

The inequality (1) leads to the following direct result:

**Assertion 2.4**. If $(X,d)$ be $p$-convex for some $p \in [1, \infty)$ then, for any $x, y, z, w \in X$, each midpoint $m(x,y) \in X$ of $x$ and $y$ and each midpoint $m(z,w) \in X$ of $z$ and $w$:

$$d(m(x,y), m(z,w)) \leq \left(\frac{1}{4}\left(d^p(x,z) + d^p(x,w) + d^p(y,z) + d^p(y,w)\right)\right)^{1/p} \tag{3}$$ □

**Definition 2.5** ($p$-Busemann curvature) [1, 2]. Suppose a metric space $(X,d)$ which admits midpoints. Then, $(X,d)$ is said to satisfy the $p$-Busemann curvature condition for some $p \in [1, \infty]$ if for each $x, y, z, w \in X$, each midpoint $m(x,y) \in X$ of $x$ and $y$ and each midpoint $m(z,w) \in X$ of $z$ and $w$, one has

$$d(m(x,y), m(z,w)) \leq \left(\frac{1}{2}\right)^{1/p} \left(d^p(x,z) + d^p(y,w)\right)^{1/p} \tag{4}$$

**Assertion 2.6**. Suppose a metric space $(X,d)$ which admits midpoints, with the midpoint map (or midset) $m: X \times X \to X$ being unique, and which satisfies the $p$-Busemann curvature condition for some $p \in [1, \infty]$. Then, one has:

$$d(m(x,y), m(z,w)) \leq \left(\frac{1}{2}\right)^{1/p} \min\left(\left(d^p(x,z) + d^p(y,w)\right)^{1/p}, \left(d^p(x,w) + d^p(y,z)\right)^{1/p}\right) \tag{5}$$

for any $x, y, z, w \in X$, where a $m(x,y) \in X$ and $m(z,w) \in X$ are, respectively, the unique midpoints of $x$ and y and $w$ and $z$. □

**Assertion 2.7**. Assume that a metric space $(X, d)$ is midpoint convex (then being 1-convex from Assertion 2.2) with unique midpoint map and that it satisfies the $p$-Busemann curvature condition for some $p \in [1, \infty]$. Then, $(X, d)$ is $p$-convex. □

The following definitions characterize near uniform convexity.

**Definition 2.8** ($\varepsilon$-separated family of points) [1]. A family of points $(x_i)_{i \in I}$ is $\varepsilon$-separated if $\inf_{i \in I} d(x_i, x_j) \geq \varepsilon$.

**Definition 2.9** (nearly uniformly convex space) [1]. A $\infty$-convex metric space $(X, d)$ is said to be nearly uniformly convex if, for any $R > 0$ and for any $\varepsilon$-separated infinite family $\{x_i\}_{i \in I}$, with $x_i \in X$, and any $y \in X$ such that $d(x_i, y) \leq r \leq R$; $\forall i \in I$, there is some $\rho = \rho(\varepsilon, R) > 0$ such that $B_{(1-\rho)r}(y) \cap c\ell(con\{x_i\}_{i \in I}) \neq \varnothing$, where $c\ell(con\{x_i\}_{i \in I})$ is the closure of the convex hull of the family $\{x_i\}_{i \in I}$.

**Definition 2.10** [5, 6]. Let $(X, d)$ be a metric space. A mapping $W: X \times X \times [0,1] \to X$ is said to be a convex structure on $X$ if for each $(x, y, t) \in X \times X \times [0,1]$ and $z \in X$,

$$d(z, W(x, y, t)) \leq t d(z, x) + (1-t) d(z, y) \tag{6}$$

**Definition 2.11** [5]. A convex metric space $(X, d, W)$ is said to be uniformly convex if, for any $\varepsilon > 0$, there exists $\delta = \delta(\varepsilon) \in (0, 1]$ such that for any $r > 0$ and $x, y, z \in X$ with $max(d(x,z), d(y,z)) \leq r$ and $d(x, y) \geq r\varepsilon$,

$$d(z, W(x, y, 1/2)) \leq r(1 - \delta) \tag{7}$$

A uniformly convex metric space $(X, d, W)$ is also referred to commonly as uniformly 1-convex, [1]. This concept may be generalized as follows:

**Definition 2.12**. A convex metric space $(X, d, W)$ is said to be uniformly $p$-convex if, for any $\varepsilon > 0$, there exists $\delta = \delta(\varepsilon) \in (0, 1]$ such that for any $r > 0$ and any $x, y, z \in X$ with $max(d(x,z), d(y,z)) \leq r$ and $d(x, y) \geq r\varepsilon$,

$$d(z, W(x, y, 1/2)) \leq r(1 - \delta)^{1/p} \tag{8}$$

**Propositions 2.13**.
1. If a convex metric space $(X, d, W)$ is uniformly convex then it is uniformly $p$-convex for any $p \geq 1$.
2. A convex metric space $(X, d, W)$ is nearly uniformly convex if, for any $\varepsilon > 0$, there exists a strictly increasing $\delta = \delta(\varepsilon) \in [0, 1]$ such that for any $r > 0$, any $z \in X$ and any $\mu = \varepsilon r$ - separated infinite family $\{x_i\}_{i \in I} \subset X$ satisfying $\sup_{i \in I} d(x_i, z) \leq r$ and

$$d(y, W(x_i, x_j, 1/2)) \leq r(1 - \delta) \; ; \quad \forall x_i, x_j (\neq x_i) \in \{x_i\}_{i \in I} \subset X \tag{9}$$

for some $y \in B_{(1-\delta)r}(z)$, where $B_\alpha(z)$ denotes an open ball of radius $\alpha$ centred at $z$.

**3**. If a convex metric space $(X, d, W)$ is uniformly convex then it is nearly uniformly convex.

**4.** If $(X,d)$ is nearly uniformly convex and strictly $\infty$-convex then for any $R > 0$ and for any $\varepsilon$-separated infinite family $\{x_i\}_{i \in I}$, with $x_i \in X$, and any $y \in X$ such that $d(x_i, y) \le r \le R$; $\forall i \in I$, there is $\rho_r = \rho_r(\varepsilon, \varepsilon_1, r) \in (0, 1 - \varepsilon_1/r]$ such that:

$d(x_i, y) \ge \varepsilon_1 = \inf_{x(\ne y) \in X} d(y, x)$, and

$c\ell \mathsf{B}_{\varepsilon_1}(y) \cap c\ell(con\{x_i\}_{i \in I}) \ne \emptyset$ and contains at most two points of $X$ such that $c\ell \mathsf{B}_{\varepsilon_1}(y) \cap c\ell(con\{x_i\}_{i \in I}) = \{w\}$ with $w \ne y$ if $y \notin \{x_i : i \in I\}$ and $c\ell \mathsf{B}_{\varepsilon_1}(y) \cap c\ell(con\{x_i\}_{i \in I}) = \{y, w\}$ if $y \in c\ell(con\{x_i\}_{i \in I})$ with a choice $\rho_r = 1 - \varepsilon_1/r$. □

## 3. Contractiveness and non-expansiveness

The following properties arise:

1) A geodesic space is a $CAT(0)$ space if and only if for any $x, y, z \in X$ and all $t \in [0, 1]$ the following inequality is satisfied:

$$d^2((1-t)x \oplus ty, z) \le (1-t)d^2(z, x) + td^2(z, y) - t(1-t)d^2(x, y) \tag{10}$$

(Proposition 1.1, [7]).

2) A $CAT(0)$ space is uniformly $p$-convex for any $p \ge 2$, [1].

3) A $CAT(0)$ space satisfies the inequalities (3) and (4) for any $p \ge 2$ since it is mid-point $p(\ge 2)$-convex for any $p \ge 1$.

4) A $CAT(0)$ space satisfies the $p$-Busemann curvature condition for any $p \ge 1$.

A general technical result involving constructions with two self-mappings in a $CAT(0)$ space follows:

**Proposition 3.1**. Let a metric space $(X, d)$ be a $CAT(0)$ space and let the mapping $Q_t : X \times X \to X$ be defined by

$$Q_t(x, y) = tx \oplus (1-t)y; \quad \forall t \in [0, 1] \tag{11}$$

for any $x, y \in X$ and let $T, S : X \to X$ be two self-mappings which satisfy the following conditions:

$$d(Tx, Ty) \le K_T d(x, y) \quad ; \quad d(Sx, Sy) \le K_S d(x, y) \tag{12}$$

for any given, some positive real constants $K_T$ and $K_S$. Then, for any given $x, y, p, q \in X$ and for any $t \in [0, 1]$, the following properties hold:

**(i)** $d^2(Q_t(Sp, Tx), Q_t(Sq, Ty)) \le (t^2 K_S^2 + (1-t)^2 K_T^2 + 4t(1-t)\min(K_S, K_T))$

$\times \max(d^2(x, y), d^2(p, q), (d(Sp, Ty), d(Sq, Ty), d(Tx, Sq), d(Sp, Tx))\min(d(p, q), d(x, y)))$ (13)

**(ii)** $d^2(Q_t(S^n p, T^n x), Q_t(S^n q, T^n y)) \le (t^2 K_S^{2n} + (1-t)^2 K_T^{2n} + 4t(1-t)\min(K_S^n, K_T^n))$

$\times \max(d^2(x, y), d^2(p, q), d(S^{n-1}p, T^{n-1}y), d(S^{n-1}q, T^{n-1}y), d(T^{n-1}x, S^{n-1}q), d(S^{n-1}p, T^{n-1}x))$

$$\times max\left(d(p,q),d(Tx,Ty)\right) \quad ; \quad \forall n \in \mathbf{Z}_+ \quad (14) \square$$

From Proposition 3.1, we get the following result:

**Theorem 3.2**. Let a metric space $(X, d)$ be a $CAT(0)$ space. Then, the following properties hold:

**(i)** Assume that $max(K_S, K_T) < 1$ (i.e. $S, T : X \to X$ are both strictly contractive). Then,

$$\lim_{n \to \infty} d\left(Q_t(S^n p, T^n x), Q_t(S^n q, T^n y)\right) = 0 \text{ for any } x, y, p, q \in X \text{ and } t \in [0, 1].$$

**(ii)** Assume that

**a)** either $K_S \in (0, 1)$ (i.e. $S : X \to X$ is strictly contractive), $K_T = 1$ (i.e. $T : X \to X$ is nonexpansive but noncontractive) and $T : X \to X$ has a fixed point, or

**b)** $K_T \in (0, 1)$ (i.e. $T : X \to X$ is strictly contractive), $K_S = 1$ (i.e. $S : X \to X$ is nonexpansive but noncontractive) and $S : X \to X$ has a fixed point.

Then,

$$\lim_{n \to \infty} \sup \left[ d^2\left(Q_t(S^n p, T^n x), Q_t(S^n q, T^n y)\right) - \min\left(t^2, (1-t)^2\right) \max\left(d^2(x, y), d^2(p, q)\right)\right] \leq 0 \quad (15)$$

for any $x, y, p, q \in X$ and $t \in [0, 1]$.

**(iii)** If $T, S : C \to C$ are both nonexpansive, where $C$ being a nonempty closed convex subset of $X$,

then $\lim_{n \to \infty} \sup \left[ d^2\left(Q_t(S^n p, T^n x), Q_t(S^n q, T^n y)\right)\right.$

$$\left. - \left(t^2 + (1-t)^2\right) \max\left(d^2(x, y), d^2(p, q)\right) + t(1-t) M(p, q, x, y) \max(d(p, q), d(x, y))\right] \leq 0 \quad (16)$$

for any $x, y, p, q \in X$ and $t \in [0, 1]$. $\square$

Proposition 3.1 leads to the following result:

**Theorem 3.5**. Let the metric space $(X, d)$ be a $CAT(0)$ space. The following properties hold:

**(i)** The following inequalities hold:

$$d\left(Q_t(S^m p, Tx), Q_t(S^m p, Ty)\right) \leq (1-t) K_T d(x, y) \; ; \forall m \in \mathbf{Z}_{0+},$$

$$d^2\left(Q_t(S^m p, T^n x), Q_t(S^m p, T^n y)\right) \leq (1-t) K_T^n d(x, y) \; ; \forall n, m \in \mathbf{Z}_{0+},$$

$$d\left(Q_t(Sp, T^m x), Q_t(Sq, T^m x)\right) \leq t K_S d(p, q) \; ; \forall m \in \mathbf{Z}_{0+}, \quad (17)$$

$$d^2\left(Q_t(S^n p, T^m x), Q_t(S^n q, T^m x)\right) \leq t K_S^n d(p, q) \; ; \forall n, m \in \mathbf{Z}_{0+}$$

**(ii)** If $(X, d)$ is complete then $Q_t : S^m p \times TC \to C$ is a strict contraction for each $t \in [0, 1]$ and each $m \in \mathbf{Z}_{0+}$, irrespective of the mapping $S : C \to C$, for any given $p, q (= p) \in X$ provided that $T : C \to C$ is strictly contractive. Thus, $\left\{d\left(Q_t(S^m p, T^n x), Q_t(S^m p, T^n y)\right)\right\} \to 0$ ; $\forall m \in \mathbf{Z}_{0+}$ as $n \to \infty$, and has a unique fixed point

$$z = z_t(m, S, p, y^*) = t S^m p \otimes (1-t) y^*$$

for each $t \in [0, 1]$ and each given $p \in C$ and $m \in \mathbf{Z}_{0+}$ and the unique fixed point $y^* = Ty^* \in C$ of $T: C \to C$ and $\{T^n y\} \to y^*$. In particular, if $t = 0$, $z_0(y^*) = y^*$ and if $t = 1$, $z_1(m, S, p, y^*) = S^m p$.

If both $S, T: C \to C$ are strictly contractive then $\{S^m p\} \to p^* = Sp^*$ with $Fix\, S = \{p^*\}$ and then
$$z = \lim_{m \to \infty} z_t(m, S, p^*, y^*) = tp^* \oplus (1-t)y^* \text{ for each } t \in [0, 1] \text{ is the unique fixed point of}$$
$Q_t: S^m p \times TC \to C$ for each $t \in [0,1]$; $\forall m \in \mathbf{Z}_{0+}$.

If $T: C \to C$ is strictly contractive and $S: C \to C$ is nonexpansive with $p^* \in Fix\, S$ ($Fix\, S = \{p^*\}$ if $S: C \to C$ is strictly contractive) then $Q_t: S^m p^* \times TC \to C$ has a unique fixed point $z_t = tp^* \oplus (1-t)y^*$ for each $t \in [0,1]$; $\forall m \in \mathbf{Z}_{0+}$.

**(iii)** If $(X, d)$ is complete then $Q_t: X \times X \to X$ is a strict contraction, irrespective of the mapping $T: C \to C$, for any given $x, y (= x) \in C$ provided that $S: C \to C$ is strictly contractive so that it has a unique fixed point and then $\{d(Q_t(S^n p, T^m x), Q_t(S^n p, T^m y))\} \to 0$; $\forall m \in \mathbf{Z}_{0+}$ as $n \to \infty$ and has a unique fixed point
$$z = z_t(m, T, p^*, x) = tp^* \oplus (1-t)T^m x$$
for each $t \in [0, 1]$ and each given $x \in C$ and $m \in \mathbf{Z}_{0+}$ and the unique existing fixed point $p^* = Sp^* \in C$ of the non-expansive mapping $S: C \to C$. □

Further useful inequalities hold as follows:
$$d^2(TQ_t(x,y), TQ_t(p,q)) = d^2(T(tx) \oplus T((1-t)y), T(tp) \oplus T((1-t)q))$$
$$d^2(Q_t(Sp,Tx), Q_t(Sq,Ty)) = d^2(tSp \oplus (1-t)Tx, tSq \oplus (1-t)Ty)$$
$$\leq td^2(Sp, tSq \oplus (1-t)Ty) + (1-t)d^2(Tx, tSq \oplus (1-t)Ty) - t(1-t)d^2(Sp, Tx)$$
$$\leq t^2 d^2(Sp, Sq) + t(1-t)d^2(Sp, Ty) - t^2(1-t)d^2(Sq, Ty)$$
$$+ t(1-t)d^2(Tx, Sq) + (1-t)^2 d^2(Tx, Ty) - t(1-t)^2 d^2(Sq, Ty)$$
$$- t(1-t)d^2(Sp, Tx)$$
$$= t^2 d^2(Sp, Sq) + (1-t)^2 d^2(Tx, Ty) + t(1-t)(d^2(Sp, Ty) + d^2(Tx, Sq))$$
$$- (t(1-t)^2 + t^2(1-t))d^2(Sq, Ty) - t(1-t)d^2(Sp, Tx)$$
$$= t^2 d^2(Sp, Sq) + (1-t)^2 d^2(Tx, Ty)$$
$$+ t(1-t)(d^2(Sp, Ty) + d^2(Tx, Sq) - d^2(Sq, Ty) - d^2(Sp, Tx))$$

$Q_t(x, y) = tx \oplus (1-t)y$; $\forall t \in [0,1]$
$$d^2(Q_t(Sp,Tx), Q_t(Sq,Ty)) = d^2(tSp \oplus (1-t)Tx, tSq \oplus (1-t)Ty)$$
$$\leq td^2(Sp, tSq \oplus (1-t)Ty) + (1-t)d^2(Tx, tSq \oplus (1-t)Ty) - t(1-t)d^2(Sp, Tx)$$

$$\leq t^2 d^2(Sp, Sq) + t(1-t)d^2(Sp, Ty) - t^2(1-t)d^2(Sq, Ty)$$
$$+ t(1-t)d^2(Tx, Sq) + (1-t)^2 d^2(Tx, Ty) - t(1-t)^2 d^2(Sq, Ty)$$
$$- t(1-t)d^2(Sp, Tx)$$
$$= t^2 d^2(Sp, Sq) + (1-t)^2 d^2(Tx, Ty) + t(1-t)\left(d^2(Sp, Ty) + d^2(Tx, Sq)\right)$$
$$- \left(t(1-t)^2 + t^2(1-t)\right)d^2(Sq, Ty) - t(1-t)d^2(Sp, Tx)$$
$$= t^2 d^2(Sp, Sq) + (1-t)^2 d^2(Tx, Ty)$$
$$+ t(1-t)\left(d^2(Sp, Ty) + d^2(Tx, Sq) - d^2(Sq, Ty) - d^2(Sp, Tx)\right)$$
$$\leq \left(t^2 K_S^2 + (1-t)^2 K_T^2\right) \max\left(d^2(p,q), d^2(x,y)\right)$$
$$+ t(1-t)\left(d(Sp, Ty) + d(Sq, Ty) + d(Tx, Sq) + d(Sp, Tx)\right)\min\left(K_S d(p,q), K_T d(x,y)\right)$$

## 4. Iterative schemes

Assume that the subsequent sequences are built:
$$x_{n+2} = Q_t(Sx_{n+1}, Tx_n) = Q_t(Sq, Ty)$$
$$x_{n+3} = Q_t(Sx_{n+2}, Tx_{n+1}) = Q_t(Sp, Tx)$$

and we assign successive values as follows:
$$x_n = y, \quad x_{n+1} = x = q, \quad x_{n+2} = p$$
so that

$$d^2(x_{n+3}, x_{n+2}) \leq \left(t^2 K_S^2 + (1-t)^2 K_T^2\right) \max\left(d^2(x_{n+1}, x_{n+2}), d^2(x_n, x_{n+1})\right)$$
$$+ t(1-t)\left(d(Sx_{n+2}, Tx_n) + d(Sx_{n+1}, Tx_n) + d(Tx_{n+1}, Sx_{n+1}) + d(Sx_{n+2}, Tx_{n+1})\right)\min\left(K_S d(x_{n+1}, x_{n+2}), K_T d(x_n, x_{n+1})\right)$$
$$\leq \left(t^2 K_S^2 + (1-t)^2 K_T^2\right) \max\left(d^2(x_{n+1}, x_{n+2}), d^2(x_n, x_{n+1})\right)$$
$$+ t(1-t)\left(d(Sx_{n+2}, Tx_n) + d(Sx_{n+1}, Tx_n) + d(Tx_{n+1}, Sx_{n+1}) + d(Sx_{n+2}, Tx_{n+1})\right)$$
$$\times \min(K_S, K_T) \min\left(d(x_{n+1}, x_{n+2}), d(x_n, x_{n+1})\right)$$
$$d^2(x_{n+3}, x_{n+2}) \leq \left(t^2 K_S^2 + (1-t)^2 K_T^2 + t(1-t)\min(K_S, K_T)\theta_n\right) \times \max\left(d^2(x_{n+1}, x_{n+2}), d^2(x_n, x_{n+1})\right)$$

The following result links an iterative scheme based on two maps to the convergence properties in CAT(0) spaces.

**Theorem 4.1**. Let the metric space $(X, d)$ be a $CAT(0)$ space and consider the iterative scheme:
$$x_{n+2} = Q_{t_n}(S_n x_{n+1}, T_n x_n) = t_n S_n x_{n+1} \oplus (1-t_n) T_n x_n; \quad \forall n \in \mathbf{Z}_{0+}$$
subject to any initial conditions $x_0, x_1 \in X$, where

$$\theta_{n+1} = \frac{d(S_n x_{n+2}, T_n x_n) + d(S_n x_{n+1}, T_n x_n) + d(T_n x_{n+1}, S_n x_{n+1}) + d(S_n x_{n+2}, T_n x_{n+1})}{\max\left(d(x_{n+1}, x_{n+2}), d(x_n, x_{n+1})\right)}; \quad \forall n \in \mathbf{Z}_{0+}$$

with $\{\theta_n\}$ and $\theta$ being depending, in general, on $x_0$ and $x_1$. The following properties hold:

**(i)** Assume that $\sup_{n \in \mathbf{Z}_{0+}} \theta_n \leq \theta < +\infty$. Thus,

a) If $\{T_n\}$ is a sequence of non-expansive (respectively, strictly contractive) sequence of mappings then there is a real sequence $\{t_n\} \subset [0,1]$ of sufficiently small elements such that $\{\rho_n\} \subset [0,1]$ (respectively, $\{\rho_n\} \subset [0,1)$),

b) If $\{S_n\}$ is a sequence of non-expansive (respectively, strictly contractive) sequence of mappings then there is a real sequence $\{t_n\} \subset [0,1]$ of elements sufficiently close to unity such that $\{\rho_n\} \subset [0,1]$ (respectively, $\{\rho_n\} \subset [0,1)$).

**(ii)** $d^2(x_{n+2}, x_{n+3}) \leq \rho_{n+1} \max\left(d^2(x_{n+1}, x_{n+2}), d^2(x_n, x_{n+1})\right); \forall n \in \mathbf{Z}_{0+}$

where

$$\rho_{n+1} = t_{n+1}^2 K_{S_{n+1}}^2 + (1-t_{n+1})^2 K_{T_{n+1}}^2 + t_{n+1}(1-t_{n+1}) \min\left(K_{S_{n+1}}, K_{T_{n+1}}\right) \theta_{n+1}; \forall n \in \mathbf{Z}_{0+}$$

for any $t_n \in [0,1]; \forall n \in \mathbf{Z}_{0+}$, and

$$d^2(x_{n+m_n}, x_{n+m_n+1}) \leq \rho_{n+m_n-1} \max\left(d^2(x_{n+m_n-1}, x_{n+m_n}), d^2(x_{n+m_n-2}, x_{n+m_n-1})\right)$$

$$\leq \left(\Pi_{i=0}^{m_n-1}[\rho_{n+i}]\right) \max\left(d^2(x_{n+1}, x_{n+2}), d^2(x_n, x_{n+1})\right); m_n \geq 2, \forall n \in \mathbf{Z}_{0+}$$

**(iii)** If $x_0, x_1 \in X$ and $x_2 = Q_{t_0}(Sx_1, Tx_0)$ are such that $d(x_1, x_2) \leq d(x_0, x_1)$ and $\rho_n \leq 1; \forall n \in \mathbf{Z}_{0+}$ then

$d(x_{n+2}, x_{n+3}) \leq d(x_n, x_{n+1}); \quad \forall n \in \mathbf{Z}_{0+}$

The above inequality is strict if $\rho_n < 1; \forall n \in \mathbf{Z}_{0+}$.